\documentclass[12pt, reqno]{amsart}
\usepackage{amsfonts}
\usepackage{amssymb}
\usepackage{amsmath}

\newtheorem{theorem}{Theorem}
\theoremstyle{plain}

\newtheorem{corollary}[theorem]{Corollary}

\newtheorem{proposition}[theorem]{Proposition}

\begin{document}
\title[] {On $p$-adic analogue of $q$-Bernstein polynomials and related integrals}
\author[]{T. Kim}
\address{Taekyun Kim. Division of General Education-Mathematics \\
Kwangwoon University, Seoul 139-701, Republic of Korea  \\}
\email{tkkim@kw.ac.kr}
\author[]{J. Choi}
\address{Jongsung Choi. Division of General Education-Mathematics \\
Kwangwoon University, Seoul 139-701, Republic of Korea  \\}
\email{jeschoi@kw.ac.kr}
\author[]{Y. H. Kim}
\address{Young-Hee Kim. Division of General
Education-Mathematics\\
Kwangwoon University, Seoul 139-701, Republic of Korea  \\}
\email{yhkim@kw.ac.kr}
\author[]{L. C. Jang}
\address{Lee-Chae Jang. Department of Mathematics and Computer Science, KonKuk University, Chungju 380-701, Republic of Korea}
\email{leechae.jang@kku.ac.kr}
\thanks{
{\it 2000 Mathematics Subject Classification}  : 11B68, 11S80,
41A30}
\thanks{\footnotesize{\it Key words and
phrases} :  $q$-Bernstein polynomial, Bernoulli numbers and
polynomials, $p$-adic $q$-integral}
\maketitle

{\footnotesize {\bf Abstract} \hspace{1mm} {Recently, T. Kim([5]) 
introduced $q$-Bernstein polynomials which are different
$q$-Bernstein polynomials of Phillips $q$-Bernstein polynomials([11, 12]). The purpose
of this paper is to study some properties of several type Kim's
$q$-Bernstein polynomials to express the $p$-adic $q$-integral of
these polynomials on $\Bbb Z_p$ associated with Carlitz's
$q$-Bernoulli numbers and polynomials. Finally, we also derive some
relations on the $p$-adic $q$-integral of the products of several
type Kim's $q$-Bernstein polynomials and the powers of them on $\Bbb
Z_p$.}

\section{Introduction}

Let $C[0, 1]$ denote the set of continuous functions on $[0, 1]$.
For $0<q<1$ and $f \in C[0, 1]$, Kim introduced the $q$-extension of
Bernstein linear operator of order $n$ for $f$ as follows:
\begin{eqnarray*}
\mathbb{B}_{n,q} (f | x)= \sum_{k=0}^{n} f(\frac{k}{n})\binom{n}{k}
[x]_q^k [1-x]_{\frac{1}{q}}^{n-k}=\sum_{k=0}^{n}
f(\frac{k}{n})B_{k,n}(x,q),
\end{eqnarray*}
where $[x]_q = \frac{1-q^x}{1-q}$, (see [5]). Here $\mathbb{B}_{n,q}
(f | x)$ is called Kim's $q$-Bernstein operator of order $n$ for
$f$. For $k, n \in \Bbb Z_+(=\Bbb N \cup \{0\})$, $B_{k,n}(x,q)=
\binom{n}{k}[x]_q^k [1-x]_{\frac{1}{q}}^{n-k}$ are called the
Kim's $q$-Bernstein polynomials of degree $n$ (see [1, 6, 11-13]).

In [2], Carlitz defined a set of numbers $\xi_k = \xi_k (q)$
inductively by
\begin{eqnarray*}
\xi_{0}=1, \quad (q \xi +1)^k -\xi_k=\left\{
\begin{array}{ll} 1 \ \ &\hbox{if}\ \ k=1,
\vspace{2mm}\\
0\ \ &\hbox{if}\ \ k>1,
\end{array}\right. \end{eqnarray*}
with the usual convention of replacing $\xi^k$ by $\xi_k$.
These numbers are $q$-analogues of ordinary Bernoulli numbers $B_k$,
but they do not remain finite for $q=1$. So he modified the definition as
follows:
\begin{eqnarray*}
\beta_{0,q}=1, \quad q(q \beta +1)^k -\beta_{k,q}=\left\{
\begin{array}{ll} 1 \ \ &\hbox{if}\ \ k=1,
\vspace{2mm}\\
0\ \ &\hbox{if}\ \ k>1,
\end{array}\right. \end{eqnarray*}
with the usual convention of replacing $\beta^k$ by $\beta_{k,q}$ (see [2]).
These numbers $\beta_{n,q}$ are called the $n$-th Carlitz
$q$-Bernoulli numbers. And Carlitz's $q$-Bernoulli
polynomials are defined by $$\beta_{k, q}(x)=(q^x \beta + [x]_q)^k =
\sum_{i=0}^k \binom{k}{i} \beta_{i, q} q^{ix} [x]_q^{k-i}.$$ As $q
\rightarrow 1$, we have $\beta_{k, q} \rightarrow B_k$ and
$\beta_{k, q}(x) \rightarrow B_k(x)$, where $B_k$ and $B_k (x)$ are
the ordinary Bernoulli numbers and polynomials, respectively.

Let $p$ be a fixed prime number. Throughout this paper, $\Bbb Z$,
$\Bbb Q$, $\Bbb Z_p$, $\Bbb Q_p$ and $\Bbb C_p$ will denote the ring
of rational integers, the field of rational numbers, the ring of
$p$-adic integers, the field of $p$-adic rational numbers and the
completion of algebraic closure of $\Bbb Q_p$, respectively.
Let $\nu_p$ be the normalized exponential valuation of $\Bbb
C_p$ such that $|p \,|_p=p^{-\nu_p (p)} =\frac{1}{p}$.

Let $q$ be regarded as either a complex number $q \in \Bbb C$ or a
$p$-adic number $q \in \Bbb C_p$. If $q\in \Bbb C$, we assume
$|q|<1$, and if $q\in \Bbb C_p$, we normally assume $|1-q|_p <1$.

We say that $f$ is a uniformly differentiable function at a point $a
\in \Bbb Z_p$ and denote this property by $f\in UD(\Bbb Z_p)$ if the difference
quotient $F_f (x,y) =\frac{f(x)-f(y)}{x-y}$ has a limit $f'(a)$ as
$(x,y)\rightarrow (a,a)$ (see [3-10]).

For $f\in UD(\Bbb Z_p)$, let us begin with the expression
\begin{eqnarray}
\frac{1}{[p^N]_q}\underset{0 \le x < p^N}{\sum}q^x f(x)= \underset{0
\le x < p^N}{\sum}f(x) \mu_q (x+p^N \Bbb Z_p),
\end{eqnarray}
representing a $q$-analogue of the Riemann sums for $f$ (see [8]).
The integral of $f$ on $\Bbb Z_p$ is defined as the limit as $n
\rightarrow \infty$ of the sums (if exists). The $p$-adic
$q$-integral on a function $f \in UD(\Bbb Z_p)$ is defined by
\begin{eqnarray*}
I_{q} (f)=\int_{\Bbb Z_p } f(x) d\mu_{q} (x) = \lim_{N \rightarrow
\infty} \frac{1}{[p^N]_q} \sum_{x=0}^{p^N-1} f(x)q^x, \quad
(\text{see [8]}).
\end{eqnarray*}

As was shown in [6], Carlitz's $q$-Bernoulli numbers can be
represented by $p$-adic $q$-integral on $\Bbb Z_p$ as follows:
\begin{eqnarray}
\int_{\Bbb Z_p } [x]_q^m d\mu_{q} (x)= \beta_{m, q}, \quad \text{for
} m \in \Bbb Z_+.
\end{eqnarray}
Also, Carlitz's $q$-Bernoulli polynomials $\beta_{k, q}(x)$ can be
represented
\begin{eqnarray}
\beta_{m, q}(x)=\int_{\Bbb Z_p } [x+y]_q^m d\mu_{q} (y), \quad
\text{for } m \in \Bbb Z_+, \quad (\text{see [6]}).
\end{eqnarray}

In this paper, we consider the $p$-adic analogue of Kim's
$q$-Bernstein polynomials on $\Bbb Z_p$ and give some properties of
the several type Kim's $q$-Bernstein polynomials to represent the
$p$-adic $q$-integral on $\Bbb Z_p$ of these polynomials. Finally,
we derive some relations on the $p$-adic $q$-integral of the
products of several type Kim's $q$-Bernstein polynomials and the
powers of them on $\Bbb Z_p$.

\medskip

\section{$q$-Bernstein polynomials associated with $p$-adic $q$-integral on $\Bbb Z_p$}

In this section, we assume that $q \in \Bbb C_p$ with $|1-q|_p<1$.

From (1), (2) and (3), we note that
\begin{eqnarray}
\int_{\Bbb Z_p}[1-x+x_1]_{\frac{1}{q}}^n d \mu_{\frac{1}{q}}(x_1)=
\frac{q^n}{(q-1)^{n-1}} \sum_{l=0}^n \binom{n}{l} (-1)^l q^{lx}
\frac{l+1}{q^{l+1}-1},\end{eqnarray}
and
\begin{eqnarray}
\int_{\Bbb Z_p}[x+x_1]_{q}^n d \mu_{q}(x_1)=
\frac{1}{(q-1)^{n-1}} \sum_{l=0}^n \binom{n}{l} (-1)^l q^{lx}
\frac{l+1}{1-q^{l+1}}.
\end{eqnarray} By (4) and (5), we get
\begin{eqnarray*}
(-1)^n q^n \int_{\Bbb Z_p} [x+x_1]_q^n d\mu_q (x_1)= \int_{\Bbb
Z_p}[1-x+x_1]_{\frac{1}{q}}^n d \mu_{\frac{1}{q}}(x_1).
\end{eqnarray*}
Therefore, we obtain the following theorem.

\begin{theorem}
For $n \in \Bbb Z_+$, we have
\begin{eqnarray*}
\int_{\Bbb Z_p}[1-x+x_1]_{\frac{1}{q}}^n d \mu_{\frac{1}{q}}(x_1)=
(-1)^n q^n \int_{\Bbb Z_p} [x+x_1]_q^n d\mu_q (x_1).
\end{eqnarray*}
\end{theorem}

\medskip

By the definition of Carlitz's $q$-Bernoulli numbers and
polynomials, we get
$$q^2 \beta_{n,q}(2)-(n+1)q^2 +q =q(q \beta+1)^n = \beta_{n,q}, \quad \text{if } n>1. $$
Thus, we have the following proposition.

\begin{proposition}
For $n \in \Bbb N$ with $n>1$, we have
\begin{eqnarray*}
\beta_{n,q}(2)=\frac{1}{q^2}\beta_{n,q}+n+1- \frac{1}{q}.
\end{eqnarray*}
\end{proposition}

\medskip
It is easy to show that
$$[1-x]_{\frac{1}{q}}^n=(1-[x]_q)^n =(-1)^n q^n [x-1]_q^n. $$
Hence, we have
$$\int_{\Bbb Z_p} [1-x]_{\frac{1}{q}}^n d \mu_{q}(x)=(-1)^n q^n \int_{\Bbb Z_p} [x-1]_q^n d \mu_q(x). $$
By (3), we get
\begin{eqnarray}
\int_{\Bbb Z_p} [1-x]_{\frac{1}{q}}^n d \mu_{q}(x)=(-1)^n q^n
\beta_{n,q}(-1).
\end{eqnarray}
By Theorem 1 and (6), we see that
\begin{eqnarray}
\int_{\Bbb Z_p} [1-x]_{\frac{1}{q}}^n d \mu_{q}(x)=(-1)^n q^n
\beta_{n,q}(-1)=\beta_{n, \frac{1}{q}}(2).
\end{eqnarray}
From (7) and Proposition 2, we have
\begin{eqnarray}
\int_{\Bbb Z_p} [1-x]_{\frac{1}{q}}^n d \mu_{q}(x)=\beta_{n,
\frac{1}{q}}(2)=q^2 \beta_{n, \frac{1}{q}}+n+1-q.
\end{eqnarray}
By (2) and (8), we obtain the following theorem.

\begin{theorem}
For $n \in \Bbb N$ with $n>1$, we have
\begin{eqnarray*}
\int_{\Bbb Z_p} [1-x]_{\frac{1}{q}}^n d \mu_{q}(x)=q^2 \int_{\Bbb
Z_p} [x]_{\frac{1}{q}}^n d \mu_{\frac{1}{q}}(x) +n+1-q.
\end{eqnarray*}
\end{theorem}

Taking the $p$-adic $q$-integral on $\Bbb Z_p$ for one Kim's
$q$-Bernstein polynomials, we get
\begin{eqnarray}
\int_{\Bbb Z_p} B_{k,n}(x,q)d \mu_{q}(x)&=& \binom{n}{k}\int_{\Bbb
Z_p}[x]_q^k
[1-x]_{\frac{1}{q}}^{n-k}d \mu_{q}(x) \\
&=&\binom{n}{k} \sum_{l=0}^{n-k} \binom{n-k}{l}(-1)^l \int_{\Bbb
Z_p}[x]_q^{k+l} d \mu_{q}(x)\notag\\
&=&\binom{n}{k} \sum_{l=0}^{n-k} \binom{n-k}{l}(-1)^l \beta_{k+l,
q}, \notag
\end{eqnarray}
and, by the $q$-symmetric property of $B_{k,n}(x,q)$, we see that
\begin{eqnarray}
\int_{\Bbb Z_p} B_{k,n}(x,q)d \mu_{q}(x)&=&\int_{\Bbb Z_p}
B_{n-k,n}(1-x, \frac{1}{q})d \mu_{q}(x)\\
&=&\binom{n}{k} \sum_{l=0}^{k} \binom{k}{l}(-1)^{k+l} \int_{\Bbb
Z_p}[1-x]_{\frac{1}{q}}^{n-l} d \mu_{q}(x). \notag
\end{eqnarray}

For $n>k+1$, by Theorem 3 and (10), we have
\begin{eqnarray}
& &\int_{\Bbb Z_p} B_{k,n}(x,q)d \mu_{q}(x) \\
& & \quad =\binom{n}{k} \sum_{l=0}^{k} (-1)^{k+l}\binom{k}{l}
[n-l+1-q+q^2 \int_{\Bbb Z_p} [x]_{\frac{1}{q}}^{n-l}d
\mu_{\frac{1}{q}}(x)]\notag\\
& & \quad =\binom{n}{k} \sum_{l=0}^{k} (-1)^{k+l}\binom{k}{l}
[n-l+1-q+q^2 \beta_{n-l, \frac{1}{q}}].\notag
\end{eqnarray}

Let $m, n, k \in \Bbb Z_+$ with $m+n>2k+1$. Then the $p$-adic
$q$-integral for the multiplication of two Kim's $q$-Bernstein
polynomials on $\Bbb Z_p$ can be given by the following relation:
\begin{eqnarray}
& & \int_{\Bbb Z_p} B_{k,n}(x,q)B_{k,m}(x,q)d \mu_{q}(x)\\& & \quad=
\binom{n}{k}\binom{m}{k}\int_{\Bbb Z_p} [x]_q^{2k}
[1-x]_{\frac{1}{q}}^{n+m-2k} d\mu_{q}(x) \notag\\
& & \quad
=\binom{n}{k}\binom{m}{k}\sum_{l=0}^{2k}\binom{2k}{l}(-1)^{l+2k}\int_{\Bbb
Z_p}[1-x]_{\frac{1}{q}}^{n+m-l} d\mu_{q}(x).\notag
\end{eqnarray}

By Theorem 3 and (12), we get
\begin{eqnarray}
& & \int_{\Bbb Z_p} B_{k,n}(x,q)B_{k,m}(x,q)d \mu_{q}(x)\\
& & \quad =\binom{n}{k}\binom{m}{k} \sum_{l=0}^{2k}
\binom{2k}{l}(-1)^{l+2k} [n+m-l+1-q+q^2 \int_{\Bbb Z_p}
[x]_{\frac{1}{q}}^{n+m-l}d \mu_{\frac{1}{q}}(x)] \notag\\
& & \quad =\binom{n}{k}\binom{m}{k} \sum_{l=0}^{2k}
\binom{2k}{l}(-1)^{l+2k} [n+m-l+1-q+q^2 \beta_{n+m-l, \frac{1}{q}}].
\notag
\end{eqnarray}

By the simple calculation, we easily get
\begin{eqnarray}
& & \int_{\Bbb Z_p} B_{k,n}(x,q)B_{k,m}(x,q)d \mu_{q}(x)\\& & \quad=
\binom{n}{k}\binom{m}{k}\int_{\Bbb Z_p} [x]_q^{2k}
[1-x]_{\frac{1}{q}}^{n+m-2k} d\mu_{q}(x)\notag \\
& & \quad
=\binom{n}{k}\binom{m}{k}\sum_{l=0}^{n+m-2k}\binom{n+m-2k}{l}(-1)^{l}\int_{\Bbb
Z_p}[x]_{q}^{l+2k} d\mu_{q}(x)\notag \\
& & \quad
=\binom{n}{k}\binom{m}{k}\sum_{l=0}^{n+m-2k}\binom{n+m-2k}{l}(-1)^{l}\beta_{l+2k,
q}.\notag
\end{eqnarray}
Continuing this process, we obtain
\begin{eqnarray}
& & \int_{\Bbb Z_p} \left(\prod_{i=1}^s B_{k,n_i}(x,q)\right)d \mu_{q}(x)  \\
& & \quad =\left(\prod_{i=1}^s \binom{n_i}{k}\right)\int_{\Bbb Z_p}
[x]_q^{sk}[1-x]_{\frac{1}{q}}^{n_1+\cdots +n_s-sk}d \mu_{q}(x) \notag\\
& & \quad =\left(\prod_{i=1}^s \binom{n_i}{k}\right)\sum_{l=0}^{sk}
\binom{s k}{l} (-1)^{sk+l}\int_{\Bbb
Z_p}[1-x]_{\frac{1}{q}}^{n_1+\cdots +n_s-l}d \mu_{q}(x).\notag
\end{eqnarray}

Let $s \in \Bbb N$ and $n_1, \ldots, n_s, k \in \Bbb Z_+$ with
$n_1+n_2+\cdots +n_s>sk+1$. By Theorem 3 and (15), we get
\begin{eqnarray}
& & \int_{\Bbb Z_p} \left(\prod_{i=1}^s B_{k,n_i}(x,q)\right)d \mu_{q}(x)  \\
& & \quad =\left(\prod_{i=1}^s \binom{n_i}{k}\right)\sum_{l=0}^{sk}
\binom{s k}{l}(-1)^{sk+l} \{\sum_{i=1}^s n_i -l+1-q+q^2 \int_{\Bbb
Z_p}[x]_{\frac{1}{q}}^{n_1+\cdots+n_s-l}d \mu_{\frac{1}{q}}(x) \}\notag \\
& & \quad =\left(\prod_{i=1}^s \binom{n_i}{k}\right)\sum_{l=0}^{sk}
\binom{s k}{l} (-1)^{sk+l}\{\sum_{i=1}^s n_i -l+1-q+q^2
\beta_{n_1+\cdots+n_s-l, \frac{1}{q}} \}. \notag
\end{eqnarray}
From the definition of binomial coefficient, we note that
\begin{eqnarray}
& & \int_{\Bbb Z_p} \left(\prod_{i=1}^s B_{k,n_i}(x,q)\right)d \mu_{q}(x) \\
& & \quad =\left(\prod_{i=1}^s \binom{n_i}{k}\right)\int_{\Bbb Z_p}
[x]_q^{sk}[1-x]_{\frac{1}{q}}^{n_1+\cdots+n_s-sk}d \mu_{q}(x) \notag \\
& & \quad =\left(\prod_{i=1}^s
\binom{n_i}{k}\right)\sum_{l=0}^{n_1+\cdots+n_s-sk}
\binom{n_1+\cdots+n_s-s k}{l}(-1)^l \int_{\Bbb Z_p}[x]_q^{sk+l}d \mu_{q}(x) \notag \\
& & \quad =\left(\prod_{i=1}^s
\binom{n_i}{k}\right)\sum_{l=0}^{n_1+\cdots+n_s-sk}
\binom{n_1+\cdots+n_s-s k}{l}(-1)^l \beta_{sk+l, q}, \notag
\end{eqnarray}
where $s \in \Bbb N$ and $n_1, \ldots, n_s, k \in \Bbb Z_+$.

By (16) and (17), we obtain the following theorem.

\begin{theorem}
(I) For $s \in \Bbb N$ and $n_1, \ldots, n_s, k \in \Bbb N$ with
$n_1+n_2+\cdots +n_s>sk+1$, we have
\begin{eqnarray*}
& & \int_{\Bbb Z_p} \left(\prod_{i=1}^s B_{k,n_i}(x,q)\right)d \mu_{q}(x) \\
& & \quad =\left(\prod_{i=1}^s \binom{n_i}{k}\right)\sum_{l=0}^{sk}
\binom{s k}{l} (-1)^{sk+l}\{\sum_{i=1}^s n_i -l+1-q+q^2
\beta_{n_1+\cdots+n_s-l, \frac{1}{q}} \}.
\end{eqnarray*}
(II) For $s \in \Bbb N$ and $n_1, \ldots, n_s, k \in \Bbb Z_+$, we
have
\begin{eqnarray*}
& & \int_{\Bbb Z_p} \left(\prod_{i=1}^s B_{k,n_i}(x,q)\right)d \mu_{q}(x) \\
& & \quad =\left(\prod_{i=1}^s
\binom{n_i}{k}\right)\sum_{l=0}^{n_1+\cdots+n_s-sk}
\binom{n_1+\cdots+n_s-s k}{l}(-1)^l \beta_{sk+l, q}. \notag
\end{eqnarray*}
\end{theorem}

\medskip

By Theorem 4, we obtain the following corollary.

\begin{corollary}
For $s \in \Bbb N$ and $n_1, \ldots, n_s, k \in \Bbb N$ with
$n_1+n_2+\cdots +n_s>sk+1$, we have
\begin{eqnarray*}
& & \sum_{l=0}^{sk} \binom{s k}{l} (-1)^{sk+l}\{\sum_{i=1}^s n_i
-l+1-q+q^2 \beta_{ n_1 +\cdots+n_s-l, \frac{1}{q}} \}
\\ & & \quad =\sum_{l=0}^{n_1+\cdots+n_s-sk} \binom{n_1+\cdots+n_s-s k}{l}(-1)^l
\beta_{sk+l, q}.
\end{eqnarray*}
\end{corollary}

\medskip

Let $s \in \Bbb N$ and $m_1, \ldots, m_s, n_1, \ldots, n_s, k \in
\Bbb Z_+$ with $m_1 n_1+ \cdots +m_s n_s>(m_1 + \cdots + m_s)k+1$.
Then we have
\begin{eqnarray}
& & \int_{\Bbb Z_p} \left(\prod_{i=1}^s B_{k,n_i}^{m_i}(x,q)\right)d \mu_{q}(x) \\
& &  =\left(\prod_{i=1}^s \binom{n_i}{k}^{m_i}\right)
\sum_{l=0}^{k\sum_{i=1}^{s}m_i}\binom{k \sum_{i=1}^{s}m_i}{l}
(-1)^{k\sum_{i=1}^{s}m_i-l} \notag \\ & & \qquad \times \int_{\Bbb
Z_p}[1-x]_{q}^{\sum_{i=1}^{s}n_i m_i -l }d \mu_{q}(x) \notag \\
& &  =\left(\prod_{i=1}^s \binom{n_i}{k}^{m_i}\right)
\sum_{l=0}^{k\sum_{i=1}^{s}m_i}\binom{k
\sum_{i=1}^{s}m_i}{l}(-1)^{k\sum_{i=1}^{s}m_i-l} \notag \\ & &
\qquad \times \{ (\sum_{i=1}^s m_i n_i-l+1)-q+q^2  \int_{\Bbb
Z_p}[x]_{\frac{1}{q}}^{\sum_{i=1}^{s}n_i
m_i  -l }d \mu_{\frac{1}{q}}(x) \} \notag \notag \\
& &  =\left(\prod_{i=1}^s \binom{n_i}{k}^{m_i}\right)
\sum_{l=0}^{k\sum_{i=1}^{s}m_i}\binom{k
\sum_{i=1}^{s}m_i}{l}(-1)^{k\sum_{i=1}^{s}m_i-l} \notag \\ & &
\qquad \times \{ (\sum_{i=1}^s m_i n_i-l+1)-q+q^2 \beta_{n_1 m_1 + n_s m_s -l,
\frac{1}{q}}\}.\notag
\end{eqnarray}
From the definition of binomial coefficient, we have
\begin{eqnarray}
& & \int_{\Bbb Z_p} \left(\prod_{i=1}^s B_{k,n_i}^{m_i}(x,q)\right)d \mu_{q}(x) \\
& &  =\left(\prod_{i=1}^s \binom{n_i}{k}^{m_i}\right)
\sum_{l=0}^{\sum_{i=1}^{s}n_i m_i-
k\sum_{i=1}^{s}m_i}\binom{\sum_{i=1}^s n_i m_i - k \sum_{i=1}^s
m_i}{l} (-1)^{l} \notag\\& & \qquad \times \int_{\Bbb
Z_p}[x]_{q}^{(m_1
+\cdots+m_s)k +l }d \mu_{q}(x) \notag \\
& & \quad =\left(\prod_{i=1}^s \binom{n_i}{k}^{m_i}\right)
\sum_{l=0}^{\sum_{i=1}^{s}n_i m_i-
k\sum_{i=1}^{s}m_i}\binom{\sum_{i=1}^s n_i m_i - k \sum_{i=1}^s
m_i}{l} \notag\\ & & \qquad \times (-1)^{l} \beta_{(m_1
+\cdots+m_s)k+l, q}. \notag
\end{eqnarray}
By (18) and (19), we obtain the following theorem.

\begin{theorem}
For $s \in \Bbb N$ and $m_1, \ldots, m_s, n_1, \ldots, n_s, k \in
\Bbb Z_+$ with $m_1 n_1+ \cdots +m_s n_s>(m_1 + \cdots + m_s)k+1$,
we have
\begin{eqnarray*}
& &\sum_{l=0}^{k\sum_{i=1}^{s}m_i}\binom{k
\sum_{i=1}^{s}m_i}{l}(-1)^{k\sum_{i=1}^{s}m_i-l}  \{ (\sum_{i=1}^s m_i
n_i-l+1)-q+q^2 \beta_{n_1 m_1 + n_s m_s -l, \frac{1}{q}}\}\\ & &
\quad=\sum_{l=0}^{\sum_{i=1}^{s}n_i m_i-
k\sum_{i=1}^{s}m_i}\binom{\sum_{i=1}^s n_i m_i - k \sum_{i=1}^s
m_i}{l} (-1)^{l} \beta_{(m_1 +\cdots+m_s)k+l, q}.
\end{eqnarray*}
\end{theorem}

\medskip

\noindent \textbf{Acknowledgement.}
This paper was supported by the research grant of Kwangwoon University in 2010.

\end{document}